\documentclass{article}
\usepackage{amsmath}
\usepackage{amsthm}
\usepackage{amssymb}
\usepackage{wasysym}
\usepackage{mathrsfs}
\usepackage{bbm} 
\usepackage{geometry} 
\usepackage[usenames, dvipsnames]{color}
\definecolor{purple}{RGB}{80, 0, 240}

\newcommand{\prob}{\Bbb P} 
\newcommand{\dpp}{d\Bbb{P}}
 
\newcommand{\E}{\Bbb E}
\DeclareMathOperator{\vol}{vol}
\DeclareMathOperator{\cone}{cone}

\DeclareMathOperator{\intt}{int}

\newcommand{\qedwhite}{\hfill \ensuremath{\Box}} 
\newcommand{\comm}[1]{} 

\newtheorem{theo}{Theorem}
\newtheorem{lemma}[theo]{Lemma}
\newtheorem{prop}[theo]{Proposition}
\newtheorem{cor}[theo]{Corollary}

\title{The Surface Area Deviation of the Euclidean Ball and a Polytope
}
\author{Steven D. Hoehner\hspace{1cm}Carsten  Sch{\"u}tt\hspace{1cm}Elisabeth M.  Werner\thanks{Partially supported by an  NSF grant}}

\date{}

\begin{document}

\maketitle
\begin{abstract}
While there is extensive literature on approximation of convex bodies by inscribed or circumscribed polytopes, much less is known in the case of  generally positioned
polytopes. Here we give upper and lower bounds for approximation of convex bodies by arbitrarily positioned polytopes  with a fixed number of vertices in the 
symmetric surface area deviation.
\end{abstract}

\section{Introduction and main results}

How well can a convex body be approximated by a
polytope?
This is a fundamental question not only  in convex geometry, but also  in view of  applications  in stochastic geometry,  
complexity,  geometric algorithms and many more (e.g., \cite{BuRe, Edelsbrunner, Gardner1995, GardnerKiderlenMilanfar, GlasauerGruber, GruberI, GruberII, PW2,  Reitzner2004II}).
\par
Accuracy of approximation is  often measured in the symmetric difference 
metric,  which reflects the volume deviation
of the approximating and approximated objects. 
Approximation of a convex body $K$ by inscribed or circumscribed polytopes 
with respect to this metric  has been studied extensively 
and many of the major questions have been resolved.  We refer to,   e.g., the surveys and  books by Gruber \cite {Gr3, Gr4, GruberBook} and the references there and to,  e.g., \cite{Ba1, Boeroetzky2000, BoeroetzkyReitzner2004, GRS1, Ludwig1999, Reitzner2004, Schneider1981, SchneiderWeil, SchuettWerner2003}.
\par
Sometimes 
it is more advantageous to consider  the surface area deviation $\Delta_s$   \cite{ BoeroetzkyCsikos2009, BoeroetzkyReitzner2004, Groemer2000} instead of the volume deviation $\Delta_v$.
It is especially desirable because
if best approximation of  convex bodies is replaced by random approximation, then
we have essentially the same amount of information for volume, surface area, and mean width (\cite{BoeroetzkyReitzner2004},\cite{BoeroetzkySchneider2010}), which are three of the quermassintegrals of a convex body (see, e.g., \cite{Gardner1995, SchneiderBook}). 
\par
For convex bodies $K$ and $L$ in $\mathbb{R}^n$ with boundaries $\partial K$ and $\partial L$,  the symmetric surface area deviation  is defined as 
\begin{equation} \label{deviation}
\Delta_s(K,L)=\vol_{n-1} \left( \partial (K \cup L)\right)-\vol_{n-1} \left( \partial (K \cap L)\right). 
\end{equation}
\par
Typically, approximation by polytopes  often involves side conditions,  like a prescribed  
number of
vertices, or, more generally,  
$k$-dimensional faces \cite{Boeroetzky2000}.   Most often in the literature,  it is required  that the body contains the approximating polytope or vice versa. This is  too restrictive as a requirement and it 
needs to be dropped. 
Here, we do exactly that  and prove upper and lower bounds for approximation of convex bodies by arbitrarily positioned polytopes in the 
symmetric surface area deviation.  This addresses  questions asked by Gruber \cite{GruberBook}.

\begin{theo}\label{upperbound:vertices}
There exists  an absolute constant $c>0$ such that for every integer $n \geq 3$, there is an $N_n$ such that for every $N \geq N_n$ there is a polytope $P_N$ in $\mathbb{R}^n$ with $N$
vertices such that
$$
\Delta_s(B_2^n,P_N)  \leq  c \,  \frac{\vol_{n-1} \left( \partial B^n_2  \right)}{N^\frac{2}{n-1}}.
$$
\end{theo}
Here, $B^n_2$ is the $n$-dimensional Euclidean unit ball with boundary $S^{n-1}=\partial B^n_2$.  
Moreover, throughout the paper $a, b, c, c_1, c_2$ will denote positive absolute constants that may change from line to line.
\par
The proof of Theorem \ref{upperbound:vertices} is based on a random construction. 
 A crucial step in its proof  is a result by J. M\"uller \cite{JMueller} on the surface deviation of a polytope {\em contained}  in the unit ball. It describes the asymptotic behavior of the surface deviation of a random polytope  $P_N$, the convex hull of $N$
 randomly (with respect to the uniform measure) and independently chosen points on the boundary of the unit ball  as the number of vertices increases. It says that 
\begin{eqnarray}\label{Mueller-surface}
& & \lim_{N\to\infty}
\frac{\vol_{n-1}(S^{n-1})-\Bbb E \vol_{n-1}(\partial P_N)}
{N^{-\frac{2}{n-1}}}   
=\frac{n-1}{n+1} \  \frac{\Gamma\left(n+\tfrac{2}{n-1}\right)}{2(n-2)!}
\frac{\left(\vol_{n-1}\left(\partial B_{2}^{n}\right)\right)^{\frac{n+1}{n-1}} }
{\left(\vol_{n-1}(B_{2}^{n-1})\right)^{\frac{2}{n-1}}}.
\end{eqnarray}

The right hand side of (\ref{Mueller-surface}) is of the order $c\,n\vol_{n-1}(\partial B^n_2)$. 
Thus, dropping the restriction that $P_N$ is contained in $B^n_2$  improves the estimate by a factor of dimension.  The same phenomenon was observed for  the volume deviation  in  \cite{LudwigSchuettWerner}. 

\vskip 3mm 
For the facets, we obtain the following lower bound for a polytope in arbitrary position.

\begin{theo}\label{AbschUntenSurf} There is a constant $c>0$ and $M_{0}\in\mathbb N$
such that for 
all $n\in\mathbb N$ with $n\geq 2$, all $M\in\mathbb N$ with $M\geq M_{0}$ and all polytopes $P_M$ in
$\mathbb R^{n}$ with no more than $M$ facets
$$
\Delta_s(B_{2}^{n},P_M)
\geq c\,\frac{\vol_{n-1}(\partial B^n_2)}{M^{\frac{2}{n-1}}}.
$$
\end{theo}
Again, we gain by a factor of dimension if we drop the requirement that the polytope contains $B^n_2$. Indeed, it follows from \cite{Gr3, McVi} that the order of best approximation  
$\Delta_v(B^n_2, P_M^{\text{best}})$ with $B^n_2 \subset P_M$ behaves asymptotically,  for $M \rightarrow \infty$,  like $\vol_{n-1}(\partial B^n_2)$. Now observe that when $B^n_2 \subset P_M$,  $n \ \Delta_v(B^n_2, P_M)= \Delta_s(B^n_2, P_M)$.
\vskip 3mm
As a corollary to Theorem \ref{AbschUntenSurf},  we deduce a lower bound in the case of simple polytopes with at most $N$ vertices.
A polytope in $\mathbb R^{n}$ is called simple if at every vertex exactly
$n$ facets meet. 
\vskip 2mm
\begin{cor}
 There is a constant $c>0$ and $N_{0}\in\mathbb N$
such that for 
all $n\in\mathbb N$ with $n\geq2$, all $N\in\mathbb N$ with $N\geq N_{0}$ and all 
simple polytopes $P_{N}$ in
$\mathbb R^{n}$ with no more than $N$ vertices
$$
\Delta_s(B_{2}^{n},P_{N})
\geq c\,\frac{\vol_{n-1}(\partial B^n_2)}{N^{\frac{2}{n-1}}}.
$$
\end{cor}

\vskip 5mm
The authors want to thank the Institute for Mathematics and Its Applications (IMA) at the University of Minnesota for their hospitality. It was during their stay there when most of the work on  the paper  was carried out. We also want to thank the referee and the editor for their careful work.

\section{Notation and auxiliary lemmas}

For a convex body $K$ in $\mathbb {R}^n$, we denote by $\intt(K)$ its interior. Its $n$-dimensional volume is $\vol_n(K)$ and the  surface area of its boundary $\partial K$ is $\vol_{n-1}(\partial K)$.
The usual surface area measure on $\partial K$ is denoted by $\mu_{\partial K}$.  The convex hull of points $x_1, \dots, x_m$ is $[x_1, \dots, x_m]$. 
\par
The affine hyperplane in $\mathbb R^{n}$ through the point $x$ and orthogonal to the vector
$\xi$ is denoted by $H(x,\xi)$.
\par
For any  further notions related to convexity, we refer to the books by e.g.,  Gruber \cite{GruberBook} and Schneider \cite{SchneiderBook}. 
\par
We start with several lemmas needed for the proof of  Theorem \ref{upperbound:vertices}. 
The first lemma says that almost all random polytopes of points chosen
from a convex body are simplicial. Intuitively this is obvious:  If we have chosen $x_{1},\dots,x_{n}$ and we want to choose $x_{n+1}$
so that it is an element of the hyperplane spanned by
$x_{1},\dots,x_{n}$, then we are choosing $x_{n+1}$  from a nullset.
We refer to, e.g., \cite{SchuettWerner2003} for the details.
\vskip 3mm
\begin{lemma}\label{Lemma: simplicial}
Almost all random polytopes of points chosen from the boundary of the
Euclidean ball with respect to the normalized surface measure are
simplicial.
\end{lemma}
\vskip 2mm
We also need the following two lemmata due to Miles \cite{Miles}.
\vskip 2mm
\begin{lemma}\label{Miles1}\rm\cite{Miles}
\begin{eqnarray*}
& & d\mu_{\partial B_{2}^{n}}(x_{1})\cdots d\mu_{\partial B_{2}^{n}}(x_{n})
\\
& & \hskip 5mm
=(n-1)!\frac{\vol_{n-1}([x_{1},\dots,x_{n}])}{(1-p^{2})^{\frac{n}{2}}}
d\mu_{\partial B_{2}^{n}\cap H}(x_{1})\cdots  d \mu_{\partial
B_{2}^{n}\cap H}(x_{n})\, dp\, d\mu_{\partial B_{2}^{n}}(\xi), 
\end{eqnarray*}
where $\xi$ is the normal to the hyperplane $H$ through $x_{1},\dots,x_{n}$ and
$p$ is the distance of the hyperplane $H$ to the origin.
\end{lemma}
\vskip 2mm

\begin{lemma}\label{Miles2}{\rm\cite{Miles}}
\begin{eqnarray*}
&&\hskip -4mm  \int_{\partial B_{2}^{n}(0,r)}\cdots\int_{\partial B_{2}^{n}(0,r)}
(\vol_n([x_{1},\dots,x_{n+1}]))^{2}   
\,d\mu_{\partial B_{2}^{n}(0,r)}(x_{1})\cdots d\mu_{\partial B_{2}^{n}
(0,r)}(x_{n+1})   \\
&& \hskip 4mm
=\frac{(n+1)r^{2n}}{n!n^{n}}
(\vol_{n-1}(\partial B_{2}^{n}(r)))^{n+1}
=\frac{(n+1)r^{n^{2}+2n-1}}{n!n^{n}}
(\vol_{n-1}(\partial B_{2}^{n}))^{n+1}.
\end{eqnarray*}
\end{lemma}
\vskip 2mm

A cap $C$ of the Euclidean ball $B_{2}^{n}$ is the intersection
of a half space $H^{-}$ with $B_{2}^{n}$. The radius of such a cap
is the radius of the $(n-1)$-dimensional ball
$B_{2}^{n}\cap H$.

\vskip 2mm
\noindent The next two ingredients needed are from \cite{SchuettWerner2003}.
\vskip 2mm
\begin{lemma}\label{SchW1}{\rm\cite{SchuettWerner2003}}
Let $H$ be a hyperplane, $p$ its distance from the origin
and $s$ the surface area of the cap
$B_{2}^{n}\cap H^{-}$, i.e., 
$$
s=\vol_{n-1}(\partial B_{2}^{n}\cap H^{-}).
$$
Then
\begin{equation*}
\frac{dp}{ds}
=-\frac{1}
{(1-p^{2})^{\frac{n-3}{2}}\vol_{n-2}(\partial B_{2}^{n-1})}.
\end{equation*}
\end{lemma}
\vskip 4mm
The following lemma is Lemma 3.13 from \cite{SchuettWerner2003}.
\vskip 2mm
\begin{lemma}\label{surface-radius}{\rm\cite{SchuettWerner2003}}
Let $C$ be a cap of the  Euclidean unit ball. Let $s$ be the
surface area of this cap and $r$ its radius. Then we have
\begin{eqnarray*}
& & \left(\frac{s}{\vol_{n-1}(B_{2}^{n-1})}\right)^{\frac{1}{n-1}}
-\frac{1}{2(n+1)}\left(\frac{s}{\vol_{n-1}(B_{2}^{n-1})}
\right)^{\frac{3}{n-1}} 
- c  \left(\frac{s}{\vol_{n-1}(B_{2}^{n-1})}
\right)^{\frac{5}{n-1}} \\
&&\\
&& \leq r(s)   \\
&& \leq\left(\frac{s}{\vol_{n-1}(B_{2}^{n-1})}\right)^{\frac{1}{n-1}}
-\frac{1}{2(n+1)}\left(\frac{s}{\vol_{n-1}(B_{2}^{n-1})}
\right)^{\frac{3}{n-1}}+ c  \left(\frac{s}{\vol_{n-1}(B_{2}^{n-1})}
\right)^{\frac{5}{n-1}}, 
\end{eqnarray*}
where $c$ is a numerical constant.
\end{lemma}
\vskip 5mm

\section{Proof of Theorem \ref{upperbound:vertices} }
To prove  Theorem \ref{upperbound:vertices}, 
we  use a probabilistic argument. We follow the strategy given in \cite{LudwigSchuettWerner}.  Instead of volume deviation,  we now have to compute the expected surface area deviation between 
$B_{2}^{n}$ and a random polytope $[x_{1},\dots,x_{N}]$
whose vertices are chosen randomly and independently from the boundary of
a Euclidean ball with slightly bigger radius. 
For technical reasons,  we choose the points from the
boundary of $B_{2}^{n}$ and we approximate $(1-\gamma)B_{2}^{n}$. 
It will turn out that  $\gamma$ is of the order $N^{-\frac{2}{n-1}}$.

\vspace{2mm}
\noindent The expected surface area difference between $(1-\gamma) B_{2}^{n}$
and a random polytope $P_{N}$ is
\begin{eqnarray*}
&&
\hskip -5mm \Bbb{E}\left[\Delta_s((1-\gamma)B_2^n,P_N)\right]  = \\
&& \hskip -5mm
 \int_{\partial B_2^n}\cdots\int_{\partial B_2^n}\biggl[\vol_{n-1}\left[\partial(P_N\cup(1-\gamma)B_2^n)\right]
-\vol_{n-1}\left[\partial(P_N\cap(1-\gamma)B_2^n)\right]\biggr]\,\dpp(x_1)\cdots\dpp(x_N),
\end{eqnarray*}
where $\prob =\frac{ \mu_{\partial B^n_2}}{\vol_{n-1}(\partial B^n_2)}$ is the uniform probability measure on $\partial B_2^n$.
For a given $N$, we choose $\gamma$ such that
\begin{equation}\label{choice c}
    \vol_{n-1}\Big(\partial\left((1-\gamma)B_2^n\right)\Big) = (1-\gamma)^{n-1}\vol_{n-1}\left(\partial B_2^n\right) 
     = \Bbb{E}\vol_{n-1}(\partial P_N).
\end{equation} 
From (\ref{Mueller-surface})  we see that for large $N$,  $(1-\gamma)^{n-1}$ is
asymptotically  equal to
\begin{equation*}\label{c-asymptotic}
1- N^{-\frac{2}{n-1}}\   \frac{n-1}{n+1}
\left(\frac{\vol_{n-1}(\partial B_{2}^{n})}
{\vol_{n-1}( B_{2}^{n-1})}
\right)^{\frac{2}{n-1}}
\frac{\Gamma(n+\frac{2}{n-1})}{2(n-2)!}.
\end{equation*}
As $(1-\gamma)^{n-1} \geq 1-(n-1) \gamma $, we get for large enough $N$ that
\begin{equation}\label{c-asymptotic2}
\gamma   \geq
  \frac{N^{-\frac{2}{n-1}}}{n+1}
\left(\frac{\vol_{n-1}(\partial B_{2}^{n})}
{\vol_{n-1}(B_{2}^{n-1})}
\right)^{\frac{2}{n-1}}
\frac{\Gamma(n+\frac{2}{n-1})}{2(n-2)!}.
\end{equation}
For $\gamma$ small enough, $(1-\gamma)^{n-1} \leq 1- (1-\frac{1}{n}) (n-1) \gamma $. Hence we get for small enough $\gamma$ and large enough $N$ that
\begin{equation}\label{c-asymptotic3}
\gamma   
\leq
 \frac{n}{n-1}   \frac{ N^{-\frac{2}{n-1}}}{n+1}
\left(\frac{\vol_{n-1}(\partial B_{2}^{n})}
{\vol_{n-1}(B_{2}^{n-1})}
\right)^{\frac{2}{n-1}}
\frac{\Gamma(n+\frac{2}{n-1})}{2(n-2)!}.
\end{equation}
Therefore, for $N$ large enough, there are absolute constants $a$ and $b$  such that
\begin{equation}\label{c-estimate}
a    \ N^{-\frac{2}{n-1}}
\leq
\gamma\leq b  \ N^{-\frac{2}{n-1}}.
\end{equation}
We continue the computation of the expected surface area deviation.  Since
$$ \vol_{n-1}\left[\partial\left((1-\gamma)B_2^n\right)\right] =\E\vol_{n-1}\left[\partial((1-\gamma)B_2^n)\cap P_N\right] + \E\vol_{n-1}\left[\partial((1-\gamma)B_2^n)\cap P_N^c\right] 
$$
and
$$\E\vol_{n-1}(\partial P_N)=\E\vol_{n-1}\left(\partial P_N\cap(1-\gamma)B_2^n\right)+\E\vol_{n-1}\left(\partial P_N\cap\left[(1-\gamma)B_2^n\right]^c\right), 
$$
our choice of $\gamma$ means that 
\begin{align} \label{choice c-2}
&\E\vol_{n-1}\left[\partial((1-\gamma)B_2^n)\cap P_N\right] + \E\vol_{n-1}\left[\partial((1-\gamma)B_2^n)\cap P_N^c\right]\\ \nonumber
=\ &\E\vol_{n-1}\left(\partial P_N\cap(1-\gamma)B_2^n\right)+\E\vol_{n-1}\left(\partial P_N\cap\left[(1-\gamma)B_2^n\right]^c\right).
\end{align} 

\noindent Thus,
\begin{align*}
&\Bbb{E}\left[\Delta_s((1-\gamma)B_2^n,P_N)\right]&\\
&=\Bbb{E}  \vol_{n-1}\left[\partial((1-\gamma)B_2^n)\cap P_N^c\right] + \Bbb{E}   \vol_{n-1}\left(\partial P_N\cap\left[(1-\gamma)B_2^n\right]^c\right)&\\
&- \Bbb{E} \vol_{n-1}\left(\partial P_N\cap(1-\gamma)B_2^n\right)- \Bbb{E} \vol_{n-1}\left[\partial((1-\gamma)B_2^n)\cap P_N\right]&\\
&=2\Big(\Bbb{E} \vol_{n-1}\left[\partial((1-\gamma)B_2^n)\cap P_N^c\right]- \Bbb{E}\vol_{n-1}\left[(1-\gamma)B_2^n\cap\partial P_N\right]\Big), &
\end{align*}
where the last equality follows from equation (\ref{choice c-2}). Hence,
\begin{align*}
&\Bbb{E}\left[\Delta_s((1-\gamma)B_2^n,P_N)\right]=&  \\
&2\int_{\partial B_2^n}\cdots\int_{\partial B_2^n}\biggl\{\vol_{n-1}\left[\partial\left((1-\gamma)B_2^n\right)\cap P_N^c\right]-\vol_{n-1}\left[(1-\gamma)B_2^n\cap\partial P_N\right] \biggr\}\,\dpp(x_1)\cdots\dpp(x_N) .&
\end{align*}
\par
\noindent We first consider
\begin{eqnarray*}
I_1&=&\int_{\partial B_2^n}\cdots\int_{\partial B_2^n}\vol_{n-1}\left[\partial\left((1-\gamma)B_2^n\right)\cap P_N^c\right]\dpp(x_1)\cdots\dpp(x_N)\\
&=& \int_{\partial B_2^n}\cdots\int_{\partial B_2^n}\vol_{n-1}\left[\partial\left((1-\gamma)B_2^n\right)\cap P_N^c\right]\mathbbm{1}_{\{0\in \intt (P_N)\}}\,\dpp(x_1)\cdots\dpp(x_N)\\
&+&\int_{\partial B_2^n}\cdots\int_{\partial B_2^n}\vol_{n-1}\left[\partial\left((1-\gamma)B_2^n\right)\cap P_N^c\right]\mathbbm{1}_{\{0\not\in \intt (P_N)\}}\,\dpp(x_1)\cdots\dpp(x_N)\\
&\leq& \int_{\partial B_2^n}\cdots\int_{\partial B_2^n}\vol_{n-1}\left[\partial\left((1-\gamma)B_2^n\right)\cap P_N^c\right]\mathbbm{1}_{\{0\in \intt (P_N)\}}\,\dpp(x_1)\cdots\dpp(x_N)\\
&+&\vol_{n-1}(\partial B_2^n)\int_{\partial B_2^n}\cdots\int_{\partial B_2^n}\mathbbm{1}_{\{0\not\in \intt (P_N)\}}\,\dpp(x_1)\cdots\dpp(x_N).
\end{eqnarray*}
\noindent By a result of \cite{Wendel} the second summand equals
$$
\vol_{n-1}\left(\partial B_{2}^{n}\right)
2^{-N+1}\sum_{k=0}^{n-1}{N-1\choose k}  
\leq\vol_{n-1}\left(\partial B_{2}^{n}\right)
2^{-N+1}n N^{n}.
$$
Therefore,
\begin{eqnarray}\label{EstI1-1}
I_{1}&\leq& \int_{\partial B_2^n}\cdots\int_{\partial B_2^n}\vol_{n-1}\left[\partial\left((1-\gamma)B_2^n\right)\cap P_N^c\right]\mathbbm{1}_{\{0\in \intt (P_N)\}}\,\dpp(x_1)\cdots\dpp(x_N)\nonumber\\
&&+\vol_{n-1}\left(\partial B_{2}^{n}\right)
2^{-N+1}n N^{n}.
\end{eqnarray}
We introduce functions $\phi_{j_1\cdots j_n}:\prod_{i=1}^N\partial B_2^n\to\Bbb{R}$ defined by 
\begin{align*}\phi_{j_1\cdots j_n}(x_1,...,x_N)=\begin{cases}
0, \text{ if }[x_{j_1},...,x_{j_n}]\text{ is not an }(n-1)\text{-dimensional face of }[x_1,...,x_N] \\
0, \text{ if }0\not\in \intt \left( [x_1,...,x_N]\right)\\
\vol_{n-1}((1-\gamma)S^{n-1}\cap P_{N}^{c}\cap \operatorname{cone}(x_{j_1},...,x_{j_n}))), \text{ otherwise}.
\end{cases}\end{align*}
\par
\noindent For vectors $y_1, \dots, y_k$ in  $\mathbb {R}^n$, 
$$
\cone(y_1, \dots, y_k) = \left\{\sum_{i=1}^k a_i y_i\bigg |\text{ }\forall i: a_i\geq 0 \right\}
$$ is the cone spanned by $y_1, \dots, y_k$.
\noindent 
From (\ref{EstI1-1}) we get
\begin{eqnarray} \label{I1,1}
I_1 &\leq& \int_{\partial B_2^n}\cdots\int_{\partial B_2^n}
\sum_{\{j_1,...,j_n\}\subset\{1,...,N\}}
\phi_{j_1,...,j_n}(x_1,...,x_N)\,\dpp(x_1)\cdots\dpp(x_N) \nonumber \\
&&+\vol_{n-1}\left(\partial B_{2}^{n}\right)
2^{-N+1}n N^{n}.
\end{eqnarray}
Inequality (\ref{I1,1}) holds since $0\in \intt (P_N)$ and
$\displaystyle\Bbb{R}^n=\bigcup_{[x_{j_1},...,x_{j_n}]\\\text{ is a facet of }P_N}\cone(x_{j_1},...,x_{j_n}).$
By Lemma \ref{Lemma: simplicial}, $P_N=[x_1,...,x_N]$ is simplicial with probability 1. Thus, the previous expression equals
\begin{align*}
{N\choose n}\int_{\partial B_2^n}\cdots\int_{\partial B_2^n}\phi_{1...n}(x_1,...,x_N)\,\dpp(x_1)\cdots\dpp(x_N) 
+\vol_{n-1}\left(\partial B_{2}^{n}\right)
2^{-N+1}n N^{n}.
\end{align*}
Let $H$ be the hyperplane containing the points $x_{1},\dots,x_{n}$.
The set of points where $H$ is not well-defined has measure $0$.
Let $H^{+}$ be the halfspace containing $0$. Then
\begin{eqnarray*}
&&\mathbb P^{N-n}
(\{(x_{n+1},\dots,x_{N})| [x_{1},\dots,x_{n}]
\hskip 1mm \mbox{is facet of}\hskip 1mm
[x_{1},\dots,x_{N}]\hskip 1mm \mbox{and}\hskip 1mm
0\in[x_{1},\dots,x_{N}]
\})  \\
&&\hskip 20mm
=\left(\frac{\mbox{vol}_{n-1}(\partial B_{2}^{n}\cap H^{+})}
{\mbox{vol}_{n-1}(\partial B_{2}^{n})}\right)^{N-n}.
\end{eqnarray*}
Therefore,  the above expression equals
\begin{align*}
&{N\choose n}\int_{\partial B_2^n}\cdots\int_{\partial B_2^n}\left[\frac{\vol_{n-1}\left(\partial B_2^n\cap H^+\right)}{\vol_{n-1}(\partial B_2^n)}\right]^{N-n}\\
&\times \vol_{n-1}\left[(1-\gamma)S^{n-1}\cap H^-\cap\cone(x_1,...,x_n)\right]\dpp(x_1)\cdots\dpp(x_n) +\vol_{n-1}\left(\partial B_{2}^{n}\right)
2^{-N+1}n N^{n}.&\\
&={N\choose n}\frac{(n-1)!}{\left(\vol_{n-1}(\partial B_2^n)\right)^n}\int_{\xi\in S^{n-1}}\int_{p=0}^1\int_{\partial(B_2^n\cap H)}\cdots\int_{\partial(B_2^n\cap H)}\left[\frac{\vol_{n-1}\left(\partial B_2^n\cap H^+\right)}{\vol_{n-1}(\partial B_2^n)}\right]^{N-n}&\\
&\times\frac{\vol_{n-1}\left([x_1,...,x_n]\right)}{\left(1-p^2\right)^{n/2}}\vol_{n-1}\left[(1-\gamma)S^{n-1}\cap H^-\cap\cone(x_1,...,x_n)\right]&\\
&\times d\mu_{\partial(B_2^n\cap H)}(x_1)\cdots d\mu_{\partial(B_2^n\cap H)}(x_n)\,dp\, d\mu_{\partial B_2^n}(\xi)
+\vol_{n-1}\left(\partial B_{2}^{n}\right)
2^{-N+1}n N^{n}.
\end{align*}
For the last equality we have used Lemma \ref{Miles1}.  It was shown in \cite{LudwigSchuettWerner}  that for $p\leq 1-\frac{1}{n}$,
\begin{eqnarray*}
\left(\frac{\vol_{n-1}(\partial
B_{2}^{n}\cap H^{+})}{\vol_{n-1}
(\partial B_{2}^{n})}\right)^{N-n} 
\leq \exp\left(-\frac{N-n}{n^{\frac{n+1}{2}}}\right)
\end{eqnarray*}
and the rest of the expression is bounded. Thus,  there is a positive constant $c_{n}$
such that for all $n\in\mathbb N$
\begin{align}\label{I1}
I_1 &\leq  {N\choose n}\frac{(n-1)!}{\left(\vol_{n-1}(\partial B_2^n)\right)^n}\int_{\xi\in S^{n-1}}\int_{p=1-\frac{1}{n} }^1\int_{\partial(B_2^n\cap H)}\cdots\int_{\partial(B_2^n\cap H)}\left[\frac{\vol_{n-1}\left(\partial B_2^n\cap H^+\right)}{\vol_{n-1}(\partial B_2^n)}\right]^{N-n}&\nonumber \\
&\times\frac{\vol_{n-1}\left([x_1,...,x_n]\right)}{\left(1-p^2\right)^{n/2}}\vol_{n-1}\left[(1-\gamma)S^{n-1}\cap H^-\cap\cone(x_1,...,x_n)\right]&\nonumber\\
&\times d\mu_{\partial(B_2^n\cap H)}(x_1)\cdots d\mu_{\partial(B_2^n\cap H)}(x_n)\, dp\, d\mu_{\partial B_2^n}(\xi)  \nonumber \\
&+\vol_{n-1}\left(\partial B_{2}^{n}\right)
2^{-N+1}n N^{n}+c_{n}\exp\left(-\frac{N-n}{n^{\frac{n+1}{2}}}\right).
\end{align}
\noindent 
Now we consider
\begin{eqnarray*}
I_2&=&\int_{\partial B_2^n}\cdots\int_{\partial B_2^n}\vol_{n-1}\left[(1-\gamma)B_2^n\cap\partial P_N\right]\dpp(x_1)\cdots\dpp(x_N)\\
&\geq& \int_{\partial B_2^n}\cdots\int_{\partial B_2^n}\vol_{n-1}\left[(1-\gamma)B_2^n\cap\partial P_N\right]\mathbbm{1}_{\{0\in \intt (P_N)\}}\,\dpp(x_1)\cdots\dpp(x_N)\\
&=&\int_{\partial B_2^n}\cdots\int_{\partial B_2^n}\sum_{\{j_1,...,j_n\}\subset\{1,...,N\}}\psi_{j_1\cdots j_n} (x_1,...,x_N)\,\dpp(x_1)\cdots\dpp(x_N), 
\end{eqnarray*}
\noindent 
where the map $\psi_{j_1\cdots j_n}:\prod_{i=1}^N\partial B_2^n\to\Bbb{R}$ is defined by \begin{align*}\psi_{j_1\cdots j_n}(x_1,...,x_N)=\begin{cases}
0, \text{ if }[x_{j_1},...,x_{j_n}]\text{ is not an }(n-1)\text{-dimensional  face of }[x_1,...,x_N] \\
0, \text{ if }0\not\in \intt \left([x_1,...,x_N]\right)\\
\vol_{n-1}\left[(1-\gamma)B_2^n\cap [x_{j_1},...,x_{j_n}]\right], \text{ otherwise}.
\end{cases}\end{align*}
\noindent 
We proceed now for $I_2$ as above for $I_1$, also using Lemma \ref{Miles1}, and get that the previous integral is greater than or equal
\begin{align} \label{I2}
&{N\choose n}\frac{(n-1)!}{\left(\vol_{n-1}(\partial B_2^n)\right)^n}\int_{\xi\in S^{n-1}}\int_{p=0}^1\int_{\partial(B_2^n\cap H)}\cdots\int_{\partial(B_2^n\cap H)}\left[\frac{\vol_{n-1}\left(\partial B_2^n\cap H^+\right)}{\vol_{n-1}(\partial B_2^n)}\right]^{N-n}& \nonumber\\
&\times\frac{\vol_{n-1}([x_1,...,x_n])}{\left(1-p^2\right)^{n/2}}\vol_{n-1}\left[(1-\gamma)B_2^n\cap H\cap\cone(x_1,...,x_n)\right]& \nonumber\\
&\times d\mu_{\partial(B_2^n\cap H)}(x_1)\cdots d\mu_{\partial(B_2^n\cap H)}(x_n)\, dp\, d\mu_{\partial B_2^n}(\xi).&
\end{align}
\noindent 
Therefore, with (\ref{I1}) and (\ref{I2}), 
\begin{align*}
&\Bbb{E}\left[\Delta_s((1-\gamma)B_2^n,P_N)\right] \leq 2{N\choose n}\frac{(n-1)!}{\left(\vol_{n-1}(\partial B_2^n)\right)^n}&\\
&\times \int_{\xi\in S^{n-1}}\int_{1-\frac{1}{n}}^1\int_{\partial(B_2^n\cap H)}\cdots\int_{\partial(B_2^n\cap H)}\left[\frac{\vol_{n-1}\left(\partial B_2^n\cap H^+\right)}{\vol_{n-1}(\partial B_2^n)}\right]^{N-n}\frac{\vol_{n-1}([x_1,...,x_n])}{\left(1-p^2\right)^{n/2}}&\\
&\times\biggl[\vol_{n-1}\left[(1-\gamma)S^{n-1}\cap H^-\cap\cone(x_1,...,x_n)\right]-\vol_{n-1}\left[(1-\gamma)B_2^n\cap [x_1,...,x_n]\right]\biggr]
&\\
&\times d\mu_{\partial(B_2^n\cap H)}(x_1)\cdots d\mu_{\partial(B_2^n\cap H)}(x_n)\, dp\, d\mu_{\partial B_2^n}(\xi)\\
&+\vol_{n-1}\left(\partial B_{2}^{n}\right)
2^{-N+1}n N^{n}+c_{n}\exp\left(-\frac{N-n}{n^{\frac{n+1}{2}}}\right).
\end{align*}
\noindent 
We notice that 
\begin{eqnarray*}
&&\hskip -10mm \vol_{n-1}\left[(1-\gamma)S^{n-1}\cap H^-\cap\cone(x_1,...,x_n)\right] \leq \\
&& \hskip 10mm \left(\frac{1-\gamma}{p}\right)^{n-1} 
\vol_{n-1}\left(  (1-\gamma)B^n_2 \cap[x_1,...,x_n] \right).
\end{eqnarray*}
Thus, 
\begin{align*}
&\Bbb{E}\left[\Delta_s((1-\gamma)B_2^n,P_N)\right]  \leq
 2{N\choose n} \frac{(n-1)!}{\left(\vol_{n-1}(\partial B_2^n)\right)^n}\\
&\times \int_{\xi\in S^{n-1}}\int_{1-\frac{1}{n}}^1\int_{\partial(B_2^n\cap H)}\cdots\int_{\partial(B_2^n\cap H)}\left[\frac{\vol_{n-1}\left(\partial B_2^n\cap H^+\right)}{\vol_{n-1}(\partial B_2^n)}\right]^{N-n} \  \frac{\left(\vol_{n-1}[x_1,...,x_n]\right)^2}{\left(1-p^2\right)^{n/2}}&\\
&\times \max\left\{0,\left(\frac{1-\gamma}{p}\right)^{n-1}-1\right\}d\mu_{\partial(B_2^n\cap H)}(x_1)\cdots d\mu_{\partial(B_2^n\cap H)}(x_n)\, dp\, d\mu_{\partial B_2^n}(\xi)
\\
&+\vol_{n-1}\left(\partial B_{2}^{n}\right)
2^{-N+1}n N^{n}+c_{n}\exp\left(-\frac{N-n}{n^{\frac{n+1}{2}}}\right).
\end{align*}
\noindent 
By Lemma \ref{Miles2} this equals
\begin{align*}
&2{N\choose n}\frac{n}{(n-1)^{n-1}}\frac{\left(\vol_{n-2}(\partial B_2^{n-1})\right)^n}{\left(\vol_{n-1}(\partial B_2^n)\right)^n}\int_{\partial B_2^n}\int_{1-\frac{1}{n}}^1\left[\frac{\vol_{n-1}\left(\partial B_2^n\cap H^+\right)}{\vol_{n-1}(\partial B_2^n)}\right]^{N-n}&\\
&\times\max\left\{0,\left(\frac{1-\gamma}{p}\right)^{n-1}-1\right\}\frac{r^{n^2-2}}{(1-p^2)^{n/2}}\,dp\, d\mu_{\partial B_2^n}(\xi)
\\
&+\vol_{n-1}\left(\partial B_{2}^{n}\right)
2^{-N+1}n N^{n}+c_{n}\exp\left(-\frac{N-n}{n^{\frac{n+1}{2}}}\right),
\end{align*}
\noindent 
where $r$ denotes the radius of $B_2^n\cap H$. The expression
$B_2^n\cap H$ is a function of the distance $p$ of the hyperplane $H$
from the origin.
Since the integral does not depend on the direction $\xi$ and $r^2+p^2=1$, this last expression is equal to 
\begin{align*}
&2{N\choose n}\frac{n}{(n-1)^{n-1}}\frac{\left(\vol_{n-2}(\partial B_2^{n-1})\right)^n}{\left(\vol_{n-1}(\partial B_2^{n-1})\right)^{n-1}}\int_{1-\frac{1}{n}}^1\left[\frac{\vol_{n-1}\left(\partial B_2^n\cap H^+\right)}{\vol_{n-1}(\partial B_2^n)}\right]^{N-n}&\\
&\times\max\left\{0,\left(\frac{1-\gamma}{p}\right)^{n-1}-1\right\}r^{n^2-n-2}\,dp
\\
&+\vol_{n-1}\left(\partial B_{2}^{n}\right)
2^{-N+1}n N^{n}+c_{n}\exp\left(-\frac{N-n}{n^{\frac{n+1}{2}}}\right), 
\end{align*}
\noindent 
which equals
\begin{align*}
&2{N\choose n}\frac{n}{(n-1)^{n-1}}\frac{\left(\vol_{n-2}(\partial B_2^{n-1})\right)^n}{\left(\vol_{n-1}(\partial B_2^n)\right)^{n-1}}\int_{1-\frac{1}{n}}^{1-\gamma}\left[1-\frac{\vol_{n-1}(\partial B_2^n\cap H^-)}{\vol_{n-1}(\partial B_2^n)}\right]^{N-n}&\\
&\times\left[\left(\frac{1-\gamma}{p}\right)^{n-1}-1\right]r^{n^2-n-2}\,dp
+\vol_{n-1}\left(\partial B_{2}^{n}\right)
2^{-N+1}n N^{n}+c_{n}\exp\left(-\frac{N-n}{n^{\frac{n+1}{2}}}\right).
\end{align*}
\noindent 
Since $p\geq 1-\frac{1}{n}$ and, by (\ref{c-estimate}), $\gamma$ is of the order $N^{-\frac{2}{n-1}}$, we have for sufficiently large $N$
$$\displaystyle\left(\frac{1-\gamma}{p}\right)^{n-1}-1\leq n(1-\gamma-p).
$$
Therefore, the previous expression can be estimated by
\begin{align*}
&2{N\choose n}\frac{n^2}{(n-1)^{n-1}}\frac{\left(\vol_{n-2}(\partial B_2^{n-1})\right)^n}{\left(\vol_{n-1}(\partial B_2^n)\right)^{n-1}}\int_{1-\frac{1}{n}}^{1-\gamma}\left[1-\frac{\vol_{n-1}(\partial B_2^n\cap H^-)}{\vol_{n-1}(\partial B_2^n)}\right]^{N-n} \frac{1-\gamma-p}{r^{n+2-n^2}}\,dp
\\
&+\vol_{n-1}\left(\partial B_{2}^{n}\right)
2^{-N+1}n N^{n}+c_{n}\exp\left(-\frac{N-n}{n^{\frac{n+1}{2}}}\right).
\end{align*}
\noindent 
Let $\phi:[0,1]\to[0,\infty)$ be the function defined by
$$
\phi(p)=\frac{\vol_{n-1}(\partial
B_{2}^{n}\cap H^{-})}{\vol_{n-1}(\partial
B_{2}^{n})}
$$
where $H$ is a hyperplane with distance $p$ from the origin.
As in \cite{LudwigSchuettWerner},  we now choose
$$
s=\phi(p)=\frac{\vol_{n-1}(\partial
B_{2}^{n}\cap H^{-})}{\vol_{n-1}(\partial
B_{2}^{n})}
$$
as our new variable under the integral.
We apply Lemma \ref{SchW1} in order to change the variable under
the integral and get that the above expression is smaller or equal to
\begin{eqnarray} \label{integral111}
&\displaystyle{N\choose n}\frac{(\vol_{n-2}(\partial
B_{2}^{n-1}))^{n-1}}
{(\vol_{n-1}(\partial
B_{2}^{n}))^{n-2}}
\frac{n^2}{(n-1)^{n-1}}
\int_{\phi(1-\gamma)}^{\phi(1-\frac{1}{n})}
(1-s)^{N-n}
(1-\gamma-p)r^{(n-1)^{2}}
ds \nonumber
\\
&+\vol_{n-1}\left(\partial B_{2}^{n}\right)
2^{-N+1}n N^{n}+c_{n}\exp\left(-\frac{N-n}{n^{\frac{n+1}{2}}}\right),
\end{eqnarray}
where $\phi(p)$ is the normalized surface area of the cap with distance $p$ 
of the hyperplane to $0$. 
Before we proceed,  we want to estimate $\phi(1-\gamma)$. 
The radius $r$ and the distance $p$ 
satisfy $1=p^{2}+r^{2}$. It was shown in \cite{LudwigSchuettWerner} that 
$$
r^{n-1}\frac{\vol_{n-1}(B_{2}^{n-1})}
{\vol_{n-1}(\partial B_{2}^{n})}
\leq \phi\left(\sqrt{1-r^{2}}\right)
\leq
\frac{1}{\sqrt{1-r^{2}}}r^{n-1}\frac{\vol_{n-1}(B_{2}^{n-1})}
{\vol_{n-1}(\partial B_{2}^{n})}.
$$
We include the argument from \cite{LudwigSchuettWerner} for completeness. We compare $\phi$ with the surface area of the
intersection $B_{2}^{n}\cap H$ of the Euclidean ball and the
hyperplane $H$. We have
$$
\frac{\vol_{n-1}(B_{2}^{n}\cap H)}
{\vol_{n-1}(\partial
B_{2}^{n})}
=r^{n-1}\frac{\vol_{n-1}(B_{2}^{n-1})}
{\vol_{n-1}(\partial
B_{2}^{n})}.
$$
Since the orthogonal projection onto $H$ maps $\partial
B_{2}^{n}\cap H^{-}$ onto $B_{2}^{n}\cap H$,  the
left hand inequality follows. 
\par
The right hand inequality follows again by
considering the orthogonal projection onto $H$. The surface area
of a
surface element of $\partial B_{2}^{n}\cap H^{-}$ equals the
surface area of the one it is mapped to in
$B_{2}^{n}\cap H$ divided by the cosine of the angle between 
the normal to $H$ and the normal to $\partial B_{2}^{n}$ at
the given point. The cosine is always greater than
$\sqrt{1-r^{2}}$.
\par
For $p=1-\gamma$ we have
$r=\sqrt{2\gamma-\gamma^{2}}\leq\sqrt{2\gamma}$.  
Therefore we get by (\ref{c-asymptotic3}), 
\begin{eqnarray}
\phi(1-\gamma)
&\leq &
\frac{2^{\frac{n-1}{2}}}{1-\gamma}
\  \frac{\vol_{n-1}(B_{2}^{n-1})}
{\vol_{n-1}(\partial B_{2}^{n})}
\left\{
\frac{n}{n-1} \frac{N^{-\frac{2}{n-1}}}{n+1}
\left(\frac{\vol_{n-1}(\partial B_{2}^{n})}
{\mbox{vol}_{n-2}(\partial B_{2}^{n-1})}
\right)^{\frac{2}{n-1}}
\frac{\Gamma(n+\frac{2}{n-1})}{2(n-2)!}\right\}^{\frac{n-1}{2}}
\nonumber\\
&=&\frac{N^{-1}}{1-\gamma}
\  \left\{\frac{n}{n+1}  \ \frac{\Gamma\left(n+\frac{2}{n-1}\right)}{(n-1)!}
\right\}^{\frac{n-1}{2}}.
\end{eqnarray} 
The quantity $\gamma$ is of the order $N^{-\frac{2}{n-1}}$, so
$1/(1-\gamma)$ is as close to $1$ as we desire for $N$ large enough.
Moreover, for all $n\in\mathbb N$
$$
\left(\frac{n}{n+1}\right)^{\frac{n-1}{2}}
\leq 1.
$$
Therefore, for all $n\in\mathbb N$ and 
$N$ large enough
$$
\phi(1-\gamma)
\leq 
\frac{1}{N}\left\{\frac{ \Gamma(n+\frac{2}{n-1})}{(n-1)!}
\right\}^{\frac{n-1}{2}}.
$$
For all $n\in\mathbb N$ with $n\geq2$, 
\begin{equation}\label{EstGamma1}
\left\{\frac{  \Gamma(n+\frac{2}{n-1})}{(n-1)!}
\right\}^{\frac{n-1}{2}}
\leq 2n.
\end{equation}
We
verify the estimate. Stirling's formula tells us that
for all $x>0$
$$
\sqrt{2\pi x}x^{x}e^{-x}
<\Gamma(x+1)
<\sqrt{2\pi x}x^{x}e^{-x}e^{\frac{1}{12x}}.
$$
Therefore,
$$
\frac{ \Gamma(n+\frac{2}{n-1})}{(n-1)!}
\leq\left(1+\frac{2}{(n-1)^2}\right)^{n-\frac{1}{2}+\frac{2}{n-1}}
\left(n-1\right)^{\frac{2}{n-1}}
e^{-\frac{2}{n-1}}e^{\frac{1}{12(n-1+\frac{2}{n-1})}}
$$
and
$$
\left(\frac{\Gamma(n+\frac{2}{n-1})}{(n-1)!}\right)^{\frac{n-1}{2}}
\leq\frac{n-1}{e}
\left(1+\frac{2}{(n-1)^2}\right)^{\frac{(n-1)(2n-1)}{4}}
\left(1+\frac{2}{(n-1)^2}\right)
e^{\frac{n-1}{24(n-1+\frac{2}{n-1})}}.
$$
The right hand expression is asymptotically equal to
$(n-1) e^{1/24}$ and (\ref{EstGamma1}) follows.
Altogether, 
\begin{equation}\label{s(1-c)}
\phi(1-\gamma)\leq \frac{2n}{N}.
\end{equation}
Since $p=\sqrt{1-r^{2}}$, 
we get for all $r$ with $0\leq r\leq 1$
$$
1-\gamma-p=1-\gamma-\sqrt{1-r^{2}}\leq \frac{1}{2}r^{2}+r^{4}-\gamma.
$$
This estimate is equivalent to 
$1-\frac{1}{2}r^{2}-r^{4}\leq\sqrt{1-r^{2}}$. The 
left hand side is negative for $r\geq 0.9$ and thus the
inequality holds for $r$ with $0.9\leq r\leq 1$. For $r$
with $0\leq r\leq 0.9$ we square both sides.
Thus the integral (\ref{integral111})  is smaller or equal to
\begin{eqnarray*}
&&{N\choose n}\frac{(\vol_{n-2}(\partial
B_{2}^{n-1}))^{n-1}}
{(\vol_{n-1}(\partial
B_{2}^{n}))^{n-2}}
\frac{n^2}{(n-1)^{n-1}} 
\int_{\phi(1-\gamma)}^{\phi(1-\frac{1}{n})}
(1-s)^{N-n}
\left(\frac{1}{2}r^{2}+r^{4}-\gamma \right)r^{(n-1)^{2}}
ds  \\
&&+\vol_{n-1}\left(\partial B_{2}^{n}\right)
2^{-N+1}n N^{n}+c_{n}\exp\left(-\frac{N-n}{n^{\frac{n+1}{2}}}\right).
\end{eqnarray*}
Now we evaluate the integral of this expression. Again, we proceed exactly as in \cite{LudwigSchuettWerner} with the obvious modifications.
We include the arguments for completeness.
We use Lemma \ref{surface-radius}. By differentiation we verify that
$(\frac{1}{2}r^{2}+r^{4}-\gamma)r^{(n-1)^{2}}$ is a monotone function
of $r$. Here we use that $\frac{1}{2}r^{2}+r^{4}-\gamma$ is
nonnegative.
\begin{eqnarray*}
&& \hskip -7mm \int_{\phi(1-\gamma)}^{\phi(1-\frac{1}{n})}
(1-s)^{N-n}
\left[\frac{1}{2}r^{2}+r^{4}-\gamma \right]r^{(n-1)^{2}}ds  
\leq 
\frac{1}{2}\int_{0}^{1}
(1-s)^{N-n}
\left[s\frac{\vol_{n-1}(\partial B_{2}^{n})}
{\vol_{n-1}(B_{2}^{n-1})}\right]^{n-1+\frac{2}{n-1}}ds  \\
&&\hskip -6mm 
+\int_{0}^{1}
(1-s)^{N-n}
\left(s\frac{\vol_{n-1}(\partial B_{2}^{n})}
{\vol_{n-1}(B_{2}^{n-1})}\right)^{n-1+\frac{4}{n-1}}ds
-\int_{0}^{1}
(1-s)^{N-n}
\gamma \left(s\frac{\vol_{n-1}(\partial B_{2}^{n})}
{\vol_{n-1}(B_{2}^{n-1})}\right)^{n-1}ds   \\ &&\hskip
10mm +\int_{0}^{\phi(1-\gamma)}
(1-s)^{N-n}
\gamma \left(s\frac{\vol_{n-1}(\partial B_{2}^{n})}
{\vol_{n-1}(B_{2}^{n-1})}\right)^{n-1}ds.
\end{eqnarray*}
By (\ref{c-asymptotic2}), 
\begin{eqnarray*}
&&\hskip -12mm\int_{\phi(1-\gamma)}^{1}
(1-s)^{N-n}
\left(\frac{1}{2}r^{2}+r^{4}-\gamma \right)r^{(n-1)^{2}}ds  \\
&&\hskip 10mm
\leq\frac{1}{2}
\left(\frac{\vol_{n-1}(\partial B_{2}^{n})}
{\vol_{n-1}(B_{2}^{n-1})}\right)^{n-1+\frac{2}{n-1}}
\frac{\Gamma(N-n+1)\Gamma(n+\frac{2}{n-1})}{\Gamma(N+1+\frac{2}{n-1})}
\nonumber\\ &&\hskip 10mm
+\left(\frac{\vol_{n-1}(\partial B_{2}^{n})}
{\vol_{n-1}(B_{2}^{n-1})}\right)^{n-1+\frac{4}{n-1}}
\frac{\Gamma(N-n+1)\Gamma(n+\frac{4}{n-1})}{\Gamma(N+1+\frac{4}{n-1})}
\nonumber\\ &&\hskip 10mm
-
\left(\frac{\vol_{n-1}(\partial
B_{2}^{n})} {\vol_{n-1}(B_{2}^{n-1})}\right)^{n-1}
\frac{\Gamma(N-n+1)\Gamma(n)}{\Gamma(N+1)}
\\ &&\hskip 25mm \times  \  
 \frac{N^{-\frac{2}{n-1}}}{n+1}
\left(\frac{\vol_{n-1}(\partial B_{2}^{n})}
{\vol_{n-1}(B_{2}^{n-1})}
\right)^{\frac{2}{n-1}}
\frac{\Gamma(n+\frac{2}{n-1})}{2(n-2)!}
\\&&\hskip 10mm
+\gamma \cdot \phi(1-\gamma )\left(\frac{\vol_{n-1}(\partial B_{2}^{n})}
{\vol_{n-1}(B_{2}^{n-1})}\right)^{n-1}
\max_{s\in[0,\phi(1-\gamma)]}(1-s)^{N-n}s^{n-1}.
\end{eqnarray*}
Thus, 
\begin{eqnarray}
&&\int_{\phi(1-\gamma)}^{1}
(1-s)^{N-n}
\left(\frac{1}{2}r^{2}+r^{4}-\gamma \right)r^{(n-1)^{2}}ds  \\
&&\hskip 10mm
\leq\frac{1}{2}
\left(\frac{\vol_{n-1}(\partial B_{2}^{n})}
{\vol_{n-1}(B_{2}^{n-1})}\right)^{n-1+\frac{2}{n-1}}
\frac{\Gamma(N-n+1)\Gamma(n+\frac{2}{n-1})}{\Gamma(N+1+\frac{2}{n-1})}
\nonumber\\ &&\hskip 10mm
+\left(\frac{\vol_{n-1}(\partial B_{2}^{n})}
{\vol_{n-1}(B_{2}^{n-1})}\right)^{n-1+\frac{4}{n-1}}
\frac{\Gamma(N-n+1)\Gamma(n+\frac{4}{n-1})}{\Gamma(N+1+\frac{4}{n-1})}
\nonumber\\ &&\hskip 10mm
-\frac{1}{2}\frac{n-1}{n+1}
\left(\frac{\vol_{n-1}(\partial
B_{2}^{n})}
{\vol_{n-1}(B_{2}^{n-1})}\right)^{n-1+\frac{2}{n-1}}
\frac{\Gamma(N-n+1)\Gamma\left(n+\frac{2}{n-1}\right)}
{\Gamma(N+1)}N^{-\frac{2}{n-1}}
\nonumber\\&&\hskip 10mm
+\gamma \cdot \phi(1-\gamma)\left(\frac{\vol_{n-1}(\partial B_{2}^{n})}
{\vol_{n-1}(B_{2}^{n-1})}\right)^{n-1}
\max_{s\in[0,\phi(1-\gamma)]}(1-s)^{N-n}s^{n-1}.
\nonumber
\end{eqnarray}
The second summand is asymptotically equal to
\begin{eqnarray}
&&\hskip -15mm \left(\frac{\vol_{n-1}(\partial B_{2}^{n})}
{\vol_{n-1}(B_{2}^{n-1})}\right)^{n-1+\frac{4}{n-1}}
\frac{(N-n)!(n-1)!n^{\frac{4}{n-1}}}{N!(N+1)^{\frac{4}{n-1}}}
\nonumber\\&& \hskip 12mm
=\left(\frac{\vol_{n-1}(\partial B_{2}^{n})}
{\vol_{n-1}(B_{2}^{n-1})}\right)^{n-1+\frac{4}{n-1}}
\frac{n^{-1+\frac{4}{n-1}}}{{N\choose n}(N+1)^{\frac{4}{n-1}}}.
\end{eqnarray}
This summand is of the order $N^{-\frac{4}{n-1}}$, 
while the others are of the order $N^{-\frac{2}{n-1}}$.
\par
We consider the sum of the first and third summands, which is equal to 
\begin{eqnarray*}
\frac{1}{2}
\left(\frac{\vol_{n-1}(\partial B_{2}^{n})}
{\vol_{n-1}(B_{2}^{n-1})}\right)^{n-1+\frac{2}{n-1}}
\frac{\Gamma(N-n+1)\Gamma(n+\frac{2}{n-1})}{\Gamma(N+1+\frac{2}{n-1})}
\left(1-
\frac{n-1}{n+1}\  
\frac{
\Gamma(N+1+\frac{2}{n-1})}
{\Gamma(N+1)N^{\frac{2}{n-1}}}\right).
\end{eqnarray*}
Since $\Gamma(N+1+\frac{2}{n-1})$ is asymptotically equal to
$(N+1)^{\frac{2}{n-1}}\Gamma(N+1)$,  the sum of the first and third
summand
is for large $N$ of the order
\begin{equation}
\frac{2}{n+1}
\left(\frac{\vol_{n-1}(\partial B_{2}^{n})}
{\vol_{n-1}(B_{2}^{n-1})}\right)^{n-1+\frac{2}{n-1}}
\frac{\Gamma(N-n+1)
\Gamma(n+\frac{2}{n-1})}{\Gamma(N+1+\frac{2}{n-1})}, 
\end{equation}
which in turn is of the order
\begin{equation}
\frac{1}{n^{2}}
\left(\frac{\vol_{n-1}(\partial B_{2}^{n})}
{\vol_{n-1}(B_{2}^{n-1})}\right)^{n-1+\frac{2}{n-1}}
{N\choose n}^{-1}N^{-\frac{2}{n-1}}.
\end{equation}
We consider now the fourth summand. By (\ref{c-estimate})
and (\ref{s(1-c)}) the fourth summand is less than
\begin{equation}\label{4summand}
bN^{-\frac{2}{n-1}}\frac{n-1}{e^{23/24}N}
\left(\frac{\vol_{n-1}(\partial B_{2}^{n})}
{\vol_{n-1}(B_{2}^{n-1})}\right)^{n-1}
\max_{s\in[0,\phi(1-\gamma)]}(1-s)^{N-n}s^{n-1}.
\end{equation}
The maximum of the function $(1-s)^{N-n}s^{n-1}$
is attained at $(n-1)/(N-1)$ and the function
is increasing on the interval $[0,(n-1)/(N-1)]$. Therefore,
since $\phi(1-\gamma)<(n-1)/(N-1)$ the
maximum of this function over the interval $[0,\phi(1-\gamma)]$ is
attained at $\phi(1-\gamma)$. By (\ref{s(1-c)}) we have
$\phi(1-\gamma)\leq e^{\frac{1}{24}}\frac{n-1}{eN}$ and thus
for $N$ sufficiently big
\begin{eqnarray*}
\max_{s\in[0,\phi(1-\gamma)]}(1-s)^{N-n}s^{n-1}
&&\leq
\left(1-\frac{n-1}{e^{23/24}N}\right)^{N-n}
\left(e^{\frac{1}{24}}\frac{n-1}{eN}\right)^{n-1}   \\
&&\leq\exp\left(\frac{n-1}{24}-\frac{(n-1)(N-n)}{e^{23/24}N}\right)
\left(\frac{n}{eN}\right)^{n-1}   \\
&&\leq\exp\left(-\frac{n-1}{4}\right)
\left(\frac{n}{eN}\right)^{n-1}.
\end{eqnarray*}
Thus we get with a new constant $b$ that (\ref{4summand}) is smaller than or equal to
$$
bN^{-\frac{2}{n-1}}\left(\frac{\vol_{n-1}(\partial B_{2}^{n})}
{\vol_{n-1}(B_{2}^{n-1})}\right)^{n-1}
e^{-\frac{n}{4}}\frac{n^{n}e^{-n}}{N^{n}}.
$$
This is asymptotically equal to
\begin{equation}
bN^{-\frac{2}{n-1}}\left(\frac{\vol_{n-1}(\partial B_{2}^{n})}
{\vol_{n-1}(B_{2}^{n-1})}\right)^{n-1}
e^{-\frac{n}{4}}\frac{1}{{N\choose n}\sqrt{2\pi n}}.
\end{equation}
Altogether, (\ref{integral111}) for $N$ sufficiently big can be estimated by
\begin{eqnarray*}
&&\hskip -5mm {N\choose n}\frac{(\vol_{n-2}(\partial
B_{2}^{n-1}))^{n-1}}
{(\vol_{n-1}(\partial
B_{2}^{n}))^{n-2}}
\frac{n^2}{(n-1)^{n-1}}
\left\{\frac{1}{n^{2}}
\left(\frac{\vol_{n-1}(\partial B_{2}^{n})}
{\vol_{n-1}(B_{2}^{n-1})}\right)^{n-1+\frac{2}{n-1}}
{N\choose n}^{-1}N^{-\frac{2}{n-1}}\right.
\\ &&\hskip 50mm\left.
+ \  bN^{-\frac{2}{n-1}}\left(\frac{\vol_{n-1}(\partial B_{2}^{n})}
{\vol_{n-1}(B_{2}^{n-1})}\right)^{n-1}
e^{-\frac{n}{4}}\frac{1}{{N\choose n}\sqrt{2\pi n}}\right\}.
\end{eqnarray*}
This can be estimated by a constant times
\begin{eqnarray} \label{last-est}
&&
(\vol_{n-1}(\partial B_{2}^{n}))
n^2\bigg\{\frac{1}{n^2}
N^{-\frac{2}{n-1}}
+bN^{-\frac{2}{n-1}}
e^{-\frac{n}{4}}\frac{1}{\sqrt{2\pi n}}\bigg\}.
\end{eqnarray}
\noindent
Finally, it should be noted that we have been estimating the approximation of $(1-\gamma) B^n_2$ and not that of $B^n_2$. Therefore we need to multiply 
(\ref{last-est}) by $(1-\gamma)^{-{(n-1)}}$. By  (\ref{c-estimate}), 
$$
(1-\gamma)^{n-1} \geq 1 - b\frac{n-1}{N^\frac{2}{n-1}},
$$
so that we have  for all $N \geq (2b(n-1))^\frac{n-1}{2}$ that $(1-\gamma)^{-{(n-1)}} \leq 2$.
\qedwhite

\bigskip

\section{Proof of  Theorem \ref{AbschUntenSurf}  }

For the proof of Theorem  \ref{AbschUntenSurf} we need several more ingredients.  Throughout this section, we denote by $\|\cdot\|_2$ the Euclidean norm on $\mathbb{R}^n$
and by $B^n_2(\xi, r)$ the $n$-dimensional  Euclidean ball with radius $r$ centered at $\xi$.
\par
For a polytope $P$, the map 
$T:\partial P \cap B^n_2 \to \partial B_{2}^{n}$ is such that it  maps  an element $x$ with a unique outer normal $N(x)$
onto the following element of $\partial B_{2}^{n}$ 
\begin{equation}\label{T} 
x \mapsto T(x)=\partial B_{2}^{n} \cap \{x+s N(x) :  s \geq 0, N(x)  \hskip 1mm \mbox{normal at}\hskip 1mm x\}.
\end{equation}
Points not having a unique normal have measure $0$ and their image is
prescribed in an arbitrary way.
\vskip 2mm
\begin{lemma}\label{AbschInnenFacet}
For all $n\in\mathbb N$ with $n\geq 2$, all $M\in\mathbb N$ with $M\geq 3$,
 all polytopes $P_M$ in $\mathbb R^{n}$ with facets $F_{i}$, $i=1,\dots,M$
and
for all $i=1,\dots,M$ we have
$$
\vol_{n-1}(T(F_{i}\cap B_{2}^{n}))
-\vol_{n-1}(F_{i}\cap  B_{2}^{n})
\geq\frac{1}{32} \  \frac{\left(\vol_{n-1}(F_{i}\cap  B_{2}^{n})\right)^{\frac{n+1}{n-1}}}{\left(\vol_{n-1}(B_{2}^{n-1})\right)^{\frac{2}{n-1}}}.
$$
\end{lemma}
\vskip 2mm
\noindent
{\bf Proof.} In the case that $F_{i}\cap B_{2}^{n}$ is the empty set, the inequality holds since 
 both sides of the inequality equal $0$.
\par
Let $\xi_{i}$, $i=1,\dots,M$, be the outer normals of $P_M$ to $F_{i}$ and let $t_{i}\in\mathbb R$
be such that $H(t_{i}\xi_{i},\xi_{i})$ is the hyperplane containing $F_{i}$.
By definition, the volume radius of $F_{i}\cap B_{2}^{n}$ is
\begin{equation}\label{AbschInnenFacet-2}
r_{i}=\left(\frac{\vol_{n-1}(F_{i}\cap  B_{2}^{n})}
{\vol_{n-1}(B_{2}^{n-1})}\right)^{\frac{1}{n-1}}.
\end{equation}
We decompose the set $F_{i}$ into the two sets
$$
F_{i}^{1}=F_{i}\cap B_{2}^{n}(t_{i}\xi_{i},\tfrac{r_{i}}{2})
\hskip 10mm\mbox{and}\hskip 10mm
F_{i}^{2}=F_{i}\cap (B_{2}^{n}(t_{i}\xi_{i},\tfrac{r_{i}}{2}))^{c}.
$$
$F_{i}^{1}$ may be the empty set but, as we shall see during the proof, 
$F_{i}^{2}$ is never empty provided $F_{i}\cap B_{2}^{n}$ is nonempty.
The map $T$ stretches an infinitesimal surface element at $x$ by the factor
$\frac{1}{|\langle \xi_{i},T(x)\rangle|}$. Therefore,
\begin{eqnarray}\label{AbschInnenFacet-3}
\vol_{n-1}(T(F_{i}\cap  B_{2}^{n}))
=\int_{F_{i}\cap  B_{2}^{n}}\frac{dx}{|\langle \xi_{i},T(x)\rangle|}
=\int_{F_{i}^{1}\cap  B_{2}^{n}}\frac{dx}{|\langle \xi_{i},T(x)\rangle|}
+\int_{F_{i}^{2}\cap  B_{2}^{n}}\frac{dx}{|\langle \xi_{i},T(x)\rangle|}.
\end{eqnarray}
For all $x\in F_{i}^{2}\cap  B_{2}^{n}$ we have 
\begin{equation}\label{AbschInnenFacet-1}
|\langle \xi_{i},T(x)\rangle|
\leq\sqrt{1-\frac{1}{4}r_{i}^{2}}.
\end{equation} 
We verify this. There is $s\geq0$ with
$T(x)=x+s\xi_{i}$. This implies $\|x+s\xi_{i}\|_{2}=1$, and consequently
$$
s+\langle x,\xi_{i}\rangle
=\sqrt{1-\|x\|_{2}^{2}+\langle x,\xi_{i}\rangle^{2}}.
$$
Moreover, $x\in (B_{2}^{n}(t_{i}\xi_{i},\tfrac{r_{i}}{2}))^{c}$ means
$$
\frac{r_{i}^{2}}{4}<\|x-t_{i}\xi_{i}\|_{2}^{2}
=\|x\|_{2}^{2}-2t_{i}\langle x,\xi_{i}\rangle+t_{i}^{2}
=\|x\|_{2}^{2}-\langle x,\xi_{i}\rangle^{2}.
$$
Thus,
$$
\langle \xi_{i},T(x)\rangle
=\langle \xi_{i},x+s\xi_{i}\rangle
=\langle \xi_{i},x\rangle+s
=\sqrt{1-\|x\|_{2}^{2}+\langle x,\xi_{i}\rangle^{2}}
<\sqrt{1-\frac{r_{i}^{2}}{4}}
$$
and we have shown (\ref{AbschInnenFacet-1}).
By (\ref{AbschInnenFacet-3}) and (\ref{AbschInnenFacet-1}),
\begin{eqnarray*}
\vol_{n-1}(T(F_{i}\cap  B_{2}^{n}))
&\geq& \vol_{n-1}(F_{i}^{1}\cap  B_{2}^{n})
+\frac{\vol_{n-1}(F_{i}^{2}\cap  B_{2}^{n})}{\sqrt{1-\frac{r_{i}^{2}}{4}}}   \\
&\geq& \vol_{n-1}(F_{i}^{1}\cap  B_{2}^{n})
+\vol_{n-1}(F_{i}^{2}\cap  B_{2}^{n})\sqrt{1+\frac{r_{i}^{2}}{4}} .
\end{eqnarray*}
Since $r_{i}\leq1$, 
\begin{eqnarray*}
\vol_{n-1}(T(F_{i}\cap  B_{2}^{n}))
&\geq&\vol_{n-1}(F_{i}^{1}\cap  B_{2}^{n})
+\left(1+\frac{r_{i}^{2}}{16} \right)\vol_{n-1}(F_{i}^{2}\cap  B_{2}^{n}) \\
&=&\vol_{n-1}(F_{i}\cap  B_{2}^{n})
+\frac{r_{i}^{2}}{16} \vol_{n-1}(F_{i}^{2}\cap  B_{2}^{n}).
\end{eqnarray*}
Since $F_{i}^{1}\subseteq B_{2}^{n}(t_{i}\xi_{i},\frac{r_{i}}{2})$,  we have
$\vol_{n-1}(F_{i}^{1})\leq\frac{r_{i}^{n-1}}{2^{n-1}}\vol_{n-1}(B_{2}^{n-1})$. 
With (\ref{AbschInnenFacet-2})
\begin{eqnarray*}
\vol_{n-1}(T(F_{i}\cap  B_{2}^{n}))
&\geq& \vol_{n-1}(F_{i}\cap  B_{2}^{n})
+\frac{r_{i}^{2}}{16} \left(\vol_{n-1}(F_{i}^{}\cap  B_{2}^{n})-\frac{r_{i}^{n-1}}{2^{n-1}}\vol_{n-1}(B_{2}^{n-1})\right) \\
&\geq&\vol_{n-1}(F_{i}\cap  B_{2}^{n})
+
\frac{\left(\vol_{n-1}(F_{i}\cap  B_{2}^{n})\right)^{\frac{n+1}{n-1}}}{ 16 \left(\vol_{n-1}(B_{2}^{n-1})\right)^{\frac{2}{n-1}} }
- 
\frac{\left(\vol_{n-1}(F_{i}\cap  B_{2}^{n})\right)^{\frac{n+1}{n-1}}}{2^{n+3}\left(\vol_{n-1}(B_{2}^{n-1})\right)^{\frac{2}{n-1}}}.
\end{eqnarray*}
Therefore,
\begin{eqnarray*}
\vol_{n-1}(T(F_{i}\cap  B_{2}^{n}))
-\vol_{n-1}(F_{i}\cap  B_{2}^{n}) 
\geq\frac{1}{32}\frac{ \left(\vol_{n-1} \left(F_{i}\cap  B_{2}^{n}\right) \right)^{\frac{n+1}{n-1}}}{\left(\vol_{n-1}(B_{2}^{n-1})\right)^{\frac{2}{n-1}}}.
\end{eqnarray*}
$\Box$
\vskip 3mm
\begin{prop}\label{AbschInnen}
For all $n\in\mathbb N$ with $n\geq 2$, all $M\in\mathbb N$ with $M\geq 3$,
 all polytopes $P_M$ in $\mathbb R^{n}$ with at most $M$ facets we have
$$
\vol_{n-1}(\partial B_{2}^{n}\cap P_{M}^{c})-
\vol_{n-1}(\partial P_{M}\cap B_{2}^{n})
\geq \frac{1}{32}
\frac{\left(\vol_{n-1}(B_{2}^{n}\cap\partial P_{M})\right)^{\frac{n+1}{n-1}}}{\left(\vol_{n-1}(B_{2}^{n-1})\right)^{\frac{2}{n-1}}} \ 
\frac{1}{M^{\frac{2}{n-1}}}.
$$
\end{prop}
\vskip 2mm
\noindent
{\bf Proof.}
Let $T$ be as in (\ref{T}). Then
\begin{eqnarray*}
\vol_{n-1}(\partial B_{2}^{n}\cap P_{M}^{c})-
\vol_{n-1}(\partial P_{M}\cap B_{2}^{n}) 
\geq\vol_{n-1}\left(\bigcup_{i=1}^{M}T(F_{i}\cap  B_{2}^{n})\right)-
\vol_{n-1}\left(\bigcup_{i=1}^{M}(F_{i}\cap B_{2}^{n})\right)   .
\end{eqnarray*}
Since the intersection of two sets $F_{i}$ and $F_{i^{\prime}}$ is a nullset
and by Lemma \ref{AbschInnenFacet}, 
\begin{eqnarray*}
&&\vol_{n-1}(\partial B_{2}^{n}\cap P_{M}^{c})-
\vol_{n-1}(\partial P_{M}\cap B_{2}^{n})  \\
&&\geq\sum_{i=1}^{M}\vol_{n-1}(T(F_{i}\cap  B_{2}^{n}))-
\sum_{i=1}^{M}\vol_{n-1}(F_{i}\cap B_{2}^{n})   
\geq\frac{1}{32}\sum_{i=1}^{M}
\frac{\left(\vol_{n-1}(F_{i}\cap  B_{2}^{n})\right)^{\frac{n+1}{n-1}}}{\left(\vol_{n-1}(B_{2}^{n-1})\right)^{\frac{2}{n-1}}}.
\end{eqnarray*}
As
$$
\sum_{i=1}^{M} \vol_{n-1}(F_{i}\cap  B_{2}^{n})=\vol_{n-1}(B_{2}^{n}\cap\partial P_{M}), 
$$
by H\"older's inequality
$$
\sum_{i=1}^{M}
\left(\vol_{n-1}(F_{i}\cap  B_{2}^{n})\right)^{\frac{n+1}{n-1}}
\geq\frac{\left(\vol_{n-1}(B_{2}^{n}\cap\partial P_{M})\right)^{\frac{n+1}{n-1}}}{M^{\frac{2}{n-1}}}.
$$
Therefore, 
$$
\vol_{n-1}(\partial B_{2}^{n}\cap P_{M}^{c})-
\vol_{n-1}(\partial P_{M}\cap B_{2}^{n}) 
\geq \frac{1}{32}
\frac{\left(\vol_{n-1}(B_{2}^{n}\cap\partial P_{M})\right)^{\frac{n+1}{n-1}}}{\left(\vol_{n-1}(B_{2}^{n-1})\right)^{\frac{2}{n-1}}}
\frac{1}{M^{\frac{2}{n-1}}}.
$$
$\Box$
\vskip 3mm
Let $R:\mathbb R^{n}\to S^{n-1}$, $ x \mapsto R(x)=\frac{x}{\|x\|_2}$ be
the radial projection. 
\vskip 2mm
\begin{lemma}\label{AbschAussen}
For all $n\in\mathbb N$ with $n\geq 2$, all $M\in\mathbb N$ with $M\geq 3$,
 all polytopes $P_{M}$ in $\mathbb R^{n}$ with $0\in \intt(P_{M})\subseteq 2B_{2}^{n}$
 and with facets $F_{i}$, $i=1,\dots,M$
and
for all $i=1,\dots,M$
$$
\vol_{n-1}(F_{i}\cap (B_{2}^{n})^{c})-
\vol_{n-1}(R(F_{i}\cap (B_{2}^{n})^{c}))
\geq \frac{1}{128}
\frac{\left(\vol_{n-1}(F_{i}\cap (B_{2}^{n})^{c})\right)^{\frac{n+1}{n-1}}}{\left(\vol_{n-1}(B_{2}^{n-1})\right)^{\frac{2}{n-1}}}.
$$
\end{lemma}
\vskip 2mm
\noindent
{\bf Proof.}
Let $\xi_{i}$, $i=1,\dots,M$, be the normals to $F_{i}$ and let $t_{i}\in\mathbb R$
be such that $H(t_{i}\xi_{i},\xi_{i})$ is the hyperplane containing $F_{i}$.
\par
Since $0$ is an interior point of $P_{M}$, $R$ maps $\partial P_{M}$ bijectively
onto $\partial B_{2}^{n}$. In particular, $R$ maps $\partial P_{M}\cap (B_{2}^{n})^{c}$
up to a nullset bijectively onto $\partial B_{2}^{n}\cap P_{M}$. The map $R$ stretches 
an infinitesimal surface element at $x$ by the factor
$\frac{\langle \xi_{i},\frac{x}{\|x\|_{2}}\rangle}{\|x\|_{2}^{n-1}}$.
\par
The volume
radius of $F_{i}\cap (B_{2}^{n})^{c}$ is
\begin{equation}\label{VR}
\rho_{i}
=\left(\frac{\vol_{n-1}(F_{i}\cap (B_{2}^{n})^{c})}{\vol_{n-1}(B_{2}^{n-1})}\right)^{\frac{1}{n-1}}.
\end{equation}
For all $x\in F_{i}\cap (B_{2}^{n})^{c}$ we have $\|x\|_{2}>1$. Thus,
\begin{equation}\label{AbschAussen-2}
\vol_{n-1}(R(F_{i}\cap (B_{2}^{n})^{c}))
=\int_{F_{i}\cap (B_{2}^{n})^{c}}
\frac{\left\langle \xi_{i},\frac{x}{\|x\|_{2}}\right\rangle}{\|x\|_{2}^{n-1}}\,dx
\leq\int_{F_{i}\cap (B_{2}^{n})^{c}}
\left\langle \xi_{i},\frac{x}{\|x\|_{2}}\right\rangle dx.
\end{equation}
We decompose the set $F_{i}\cap (B_{2}^{n})^{c}$
into two sets
$$
A_{i}=F_{i}\cap (B_{2}^{n})^{c}\cap B_{2}^{n}\left(t_{i}\xi_{i},\tfrac{\rho_{i}}{2}\right)
\hspace{5mm} \text{and} \hspace{5mm}
B_{i}=F_{i}\cap (B_{2}^{n})^{c}\cap \left(B_{2}^{n}\left(t_{i}\xi_{i},\tfrac{\rho_{i}}{2}\right)\right)^{c}.
$$
For all $x\in F_{i}\cap (B_{2}^{n}(t_{i}\xi_{i},\frac{\rho_{i}}{2}))^{c}$ we have 
\begin{equation}\label{AbschAussen-1}
\left\langle \xi_{i},\frac{x}{\|x\|_{2}}\right\rangle
\leq\sqrt{1-\frac{\rho_{i}^{2}}{4\|x\|_{2}^{2}}}.
\end{equation}
We verify this. The inequality $\|x-t_{i}\xi_{i}\|_{2}>\frac{\rho_{i}}{2}$ implies
$$
\frac{\rho_{i}^{2}}{4}<\|x\|_{2}^{2}-2t_{i}\langle x,\xi_{i}\rangle+t_{i}^{2}
=\|x\|_{2}^{2}-\langle x,\xi_{i}\rangle^{2}.
$$
Thus (\ref{AbschAussen-1}) follows.
By (\ref{AbschAussen-2}) and (\ref{AbschAussen-1}),
$$
\vol_{n-1}(R(F_{i}\cap (B_{2}^{n})^{c}))
\leq\int_{A_{i}}
\left\langle \xi_{i},\frac{x}{\|x\|_{2}}\right\rangle dx
+\int_{B_{i}}
\left\langle \xi_{i},\frac{x}{\|x\|_{2}}\right\rangle dx  
\leq\int_{A_{i}} dx
+\int_{B_{i}}
\sqrt{1-\frac{\rho_{i}^{2}}{4\|x\|_{2}^{2}}}\, dx   .
$$
Since $P_{M}\subseteq 2B_{2}^{n}$,
\begin{eqnarray*}
\vol_{n-1}(R(F_{i}\cap (B_{2}^{n})^{c}))
&\leq& \vol_{n-1}(A_{i})+\vol_{n-1}(B_{i})\sqrt{1-\frac{\rho_{i}^{2}}{16}}  \\
&\leq&\vol_{n-1}(F_{i}\cap (B_{2}^{n})^{c})-\frac{\rho_{i}^{2}}{64} \vol_{n-1}(B_{i}).
\end{eqnarray*}
Since $\vol_{n-1}(A_{i})\leq\frac{\rho_{i}^{n-1}}{2^{n-1}} \vol_{n-1}(B_{2}^{n-1})$, 
we have $ \vol_{n-1}(B_{i}) \geq  \vol_{n-1}(F_{i}\cap (B_{2}^{n})^{c}) -\frac{\rho_{i}^{n-1}}{2^{n-1}} \vol_{n-1}(B_{2}^{n-1})$.
Therefore, with (\ref{VR}), 
\begin{eqnarray*}
\vol_{n-1}(R(F_{i}\cap (B_{2}^{n})^{c}))
&\leq&\left(1-\frac{\rho_{i}^{2}}{64}\right) \vol_{n-1}(F_{i}\cap (B_{2}^{n})^{c})
+\frac{\rho_{i}^{n+1}}{2^{n+5}} \vol_{n-1}(B_{2}^{n-1})  \\
&=& \vol_{n-1}(F_{i}\cap (B_{2}^{n})^{c})
-\frac{\left(\vol_{n-1}(F_{i}\cap (B_{2}^{n})^{c})\right)^{\frac{n+1}{n-1}}}{\left(\vol_{n-1}(B_{2}^{n-1})\right)^{\frac{2}{n-1}}}
\left(\frac{1}{64}-\frac{1}{2^{n+5}}\right).
\end{eqnarray*}
$\Box$
\vskip 3mm
\begin{prop}\label{AbschAuss} 
For all $n\in\mathbb N$ with $n\geq 2$, all $M\in\mathbb N$ with $M\geq 3$,
 all polytopes $P_{M}$ in $\mathbb R^{n}$ with at most $M$ facets and with
$0\in \intt(P_{M})\subseteq 2B_{2}^{n}$
$$
\vol_{n-1}(\partial  P_{M}\cap (B_{2}^{n})^{c})-
\vol_{n-1}(\partial B_{2}^{n}\cap P_{M})
\geq \frac{1}{128}
\frac{\left(\vol_{n-1}(\partial P_{M}\cap (B_{2}^{n})^{c})\right)^{\frac{n+1}{n-1}}}
{\left(\vol_{n-1}(B_{2}^{n-1})\right)^{\frac{2}{n-1}}M^{\frac{2}{n-1}}}.
$$
\end{prop}
\vskip 2mm
\noindent
{\bf Proof.} By Lemma \ref{AbschAussen}, 
\begin{eqnarray*}
&&\vol_{n-1}(\partial  P_{M}\cap (B_{2}^{n})^{c})-
\vol_{n-1}(\partial B_{2}^{n}\cap P_{M})   \\
&&\geq\sum_{i=1}^{M}\bigg[\vol_{n-1}(F_{i}\cap (B_{2}^{n})^{c})-
\vol_{n-1}(R(F_{i}\cap (B_{2}^{n})^{c}))\bigg]  
\geq\frac{1}{128}\sum_{i=1}^{M}
\frac{\left(\vol_{n-1}(F_{i}\cap (B_{2}^{n})^{c})\right)^{\frac{n+1}{n-1}}}{\left(\vol_{n-1}(B_{2}^{n-1})\right)^{\frac{2}{n-1}}}.
\end{eqnarray*}
As
$$
\vol_{n-1}(\partial P_{M}\cap (B_{2}^{n})^{c})
=\sum_{i=1}^{M} \vol_{n-1}(F_{i}\cap (B_{2}^{n})^{c}), 
$$
 H\"older's inequality implies
$$
\left(\sum_{i=1}^{M}
\left(\vol_{n-1}(F_{i}\cap (B_{2}^{n})^{c})\right)^{\frac{n+1}{n-1}}\right)^{\frac{n-1}{n+1}}
M^{\frac{2}{n+1}}\geq\sum_{i=1}^{M}
\vol_{n-1}(F_{i}\cap (B_{2}^{n})^{c})
=\vol_{n-1}(\partial P_{M}\cap (B_{2}^{n})^{c}).
$$
Consequently,
\begin{eqnarray*}
\vol_{n-1}(\partial  P_{M}\cap (B_{2}^{n})^{c})-
\vol_{n-1}(\partial B_{2}^{n}\cap P_{M})   
\geq\frac{1}{128}
\frac{\left(\vol_{n-1}(\partial P_{M}\cap (B_{2}^{n})^{c})\right)^{\frac{n+1}{n-1}}}
{\left(\vol_{n-1}(B_{2}^{n-1})\right)^{\frac{2}{n-1}}M^{\frac{2}{n-1}}}.
\end{eqnarray*}
$\Box$
\vskip 3mm
\noindent
{\bf Proof of Theorem \ref{AbschUntenSurf}.}
We may assume that the origin is an interior point of $P_{M}$.  If not, then
$P_{M}$ is contained in a Euclidean half ball and, by convexity,
the surface area of $P_{M}$ is smaller than that of the half ball, 
$ \vol_{n-1}\left(\partial P_{M} \right) \leq\frac{1}{2} \vol_{n-1}(\partial B_{2}^{n})+\vol_{n-1}(B_{2}^{n-1})$. So,
for sufficiently large $M$, 
\begin{eqnarray*}
\Delta_s(B_{2}^{n},P_{M})
\geq \vol_{n-1}(\partial B^n_2)
- \vol_{n-1}(\partial P_{M})
\geq\tfrac{1}{2}\vol_{n-1}(\partial B^n_2)- \vol_{n-1}(B_{2}^{n-1})
\geq \frac{\vol_{n-1}(\partial B^n_2)} {M^{\frac{2}{n-1}}}.
\end{eqnarray*}
In the same way, we see that for sufficiently large $M$ we may assume that
$\vol_{n-1}(\partial P_{M})\geq\frac{1}{2} \vol_{n-1}(\partial B_{2}^{n})$.
\par
Moreover, we may assume that $P_{M}\subseteq 2B_{2}^{n}$.
If not, there is $x_{0}\in P_{M}$ with $\|x_{0}\|_{2}\geq 2$. 
For $M$ sufficiently big we may assume that 
$\frac{1}{2}B_{2}^{n}\subseteq P_{M}$. Therefore,
$$
\Delta_s(B_{2}^{n},P_{M})
\geq \vol_{n-1}(\partial[x_{0},\tfrac{1}{2}B_{2}^{n}]\cap (B_{2}^{n})^{c}),
$$
where $[x_{0},\tfrac{1}{2}B_{2}^{n}]$ denotes the convex hull of the point $x_{0}$
with the Euclidean ball of radius $\frac{1}{2}$.
\par
By Propositions \ref{AbschInnen} and \ref{AbschAuss}, 
\begin{eqnarray*}
\Delta_s(B_{2}^{n},P_{M})
&=&\vol_{n-1}(\partial B_{2}^{n} \cap P_{M}^{c})
- \vol_{n-1}(\partial P_{M}\cap B_{2}^{n})     \\
&&+\vol_{n-1}(\partial P_{M}\cap (B_{2}^{n})^{c})
-\vol_{n-1}(\partial B_{2}^{n}\cap P_{M})   \\
&\geq&\frac{1}{32} \  
\frac{\left(\vol_{n-1}(B_{2}^{n}\cap\partial P_{M})\right) ^{\frac{n+1}{n-1}}} {\left(\vol_{n-1}(B_{2}^{n-1})\right)^{\frac{2}{n-1}}M^{\frac{2}{n-1}}}
+\frac{1}{128} \  \frac{\left(\vol_{n-1}(\partial P_{M}\cap (B_{2}^{n})^{c})\right)^{\frac{n+1}{n-1}}}
{\left(\vol_{n-1}(B_{2}^{n-1})\right)^{\frac{2}{n-1}}M^{\frac{2}{n-1}}}     .
\end{eqnarray*}
By H\"older's inequality, 
\begin{eqnarray*}
\Delta_s(B_{2}^{n},P_{M})
\geq\frac{1}{128 \cdot2^{\frac{2}{n-1}}} \  
\frac{\left(\vol_{n-1}(\partial P_{M})\right)^{\frac{n+1}{n-1}}} {\left(\vol_{n-1}(B_{2}^{n-1})\right)^{\frac{2}{n-1}}M^{\frac{2}{n-1}}}.
\end{eqnarray*}
For sufficiently large $M$ we have $\vol_{n-1}(\partial P_{M})\geq\frac{1}{2}\vol_{n-1}(\partial B^n_2)$.
Therefore,
$$
\Delta_s(B_{2}^{n},P_{M})\geq
\frac{1}{2^{12}} \  
\frac{\left(\vol_{n-1}(\partial B^n_2) \right)^{\frac{n+1}{n-1}}} {\left(\vol_{n-1}(B_{2}^{n-1})\right)^{\frac{2}{n-1}}M^{\frac{2}{n-1}}}.
$$
There is a constant $c>0$ such that for all $n\in\mathbb N$ with $n\geq 2$
$$
\left(\frac{\vol_{n-1}(\partial B^n_2)}{\vol_{n-1}(B_{2}^{n-1})}\right)^{\frac{2}{n-1}}\geq c.
$$
Therefore, with a new constant $c$, 
$$
\Delta_s(B_{2}^{n},P_{M})
\geq
c\,\frac{\vol_{n-1}(\partial B^n_2)}{M^{\frac{2}{n-1}}}.
$$
$\Box$

\noindent 
Steven Hoehner\\
{\small Department of Mathematics } \\
{\small Case Western Reserve University } \\
{\small Cleveland, Ohio 44106, U. S. A. } \\
{\small \tt sdh60@case.edu}\\ \\
\and
Carsten Sch\"utt\\
{\small Mathematisches Institut}\\
{\small Universit\"at Kiel}\\
{\small 24105 Kiel, Germany}\\
{\small \tt schuett@math.uni-kiel.de }\\ \\
\noindent
\and 
Elisabeth Werner\\
{\small Department of Mathematics \ \ \ \ \ \ \ \ \ \ \ \ \ \ \ \ \ \ \ Universit\'{e} de Lille 1}\\
{\small Case Western Reserve University \ \ \ \ \ \ \ \ \ \ \ \ \ UFR de Math\'{e}matique }\\
{\small Cleveland, Ohio 44106, U. S. A. \ \ \ \ \ \ \ \ \ \ \ \ \ \ 59655 Villeneuve d'Ascq, France}\\
{\small \tt elisabeth.werner@case.edu}\\ \\

\end{document}